\documentclass[a4paper,12pt]{article}
\usepackage[utf8]{inputenc}
\usepackage{amsmath}  
\usepackage{amssymb}       
\usepackage{amsfonts}
\usepackage{dsfont}
\usepackage{times}
\usepackage{amsthm} 
\usepackage{epsfig}

\usepackage{xcolor}

\newcommand{\DATUM}{06.04.2020}       
\newcommand{\one}{\mathbf{1}}

\newcommand{\cD}{\mathcal{D}}

\newcommand{\cF}{\mathcal{F}}
\newcommand{\cG}{\mathcal{G}}

\newcommand{\field}[1]{\mathbb{#1}}
 
\newcommand{\RR}{\field{R}}     
\newcommand{\NN}{\field{N}}     
\newcommand{\cirS}{\mathop{\bigcirc\kern -.73em {\scriptstyle{\rm S}}}}

\newtheorem{theorem}{Theorem}
\newtheorem{lemma}[theorem]{Lemma}             
\newtheorem{definition}[theorem]{Definition}   
\newtheorem{proposition}[theorem]{Proposition} 
\theoremstyle{plain}
\begin{document}
\bibliographystyle{plain}
\title{A Variational Approach to a $L^1$-Minimization Problem Based on the Milman-Pettis Theorem}

\author{Alexander~Hach $<$a.hach@tu-bs.de$>$, \\[1ex]
TU Braunschweig \\ Universitätsplatz~2, 38106 Braunschweig, Germany}

\date{\DATUM}

\maketitle

\begin{abstract}
\noindent We develop a variational approach to the minimization problem of functionals of the type $\frac12\left\lVert \nabla \phi \right\rVert^2_2 + \beta \left\lVert \phi \right\rVert_1$ constrained by $\left\lVert \phi \right\rVert_2 = 1$ which is related to the characterization of cases satisfying the sharp Nash inequality. Employing theory of uniform convex spaces by Clarkson and the Milman-Pettis theorem we are able account for the non-reflexivity of $L^1$ and implement the direct method of calculus of variations. By deriving the Euler-Lagrange equation we verify that the minimizers are up to rearrangement compactly supported solutions to the inhomogeneous Helmholtz equation and we study their scaling behaviour in $\beta$. 
\end{abstract}

\smallskip
\noindent \textbf{Keywords}:  calculus of variations $\cdot$ $L^1$-minimization $\cdot$ uniform convexity $\cdot$ Milman-Pettis theorem $\cdot$ Nash inequality
\bigskip
\thispagestyle{empty}

\setcounter{page}{1}


We consider the minimization problem for a family $\left(\cF_\beta\right)_\beta$ of linear functionals defined on the intersection of Sobolev space $H^1(\RR^3) := W^{1,2}(\RR^3)$ and Lebesgue space $L^1(\RR^3)$ given by $\cF_\beta(\phi)  = \frac{1}{2}\left\lVert\nabla\phi\right\rVert_2^2 + \beta \left\lVert \phi \right\rVert_1$ under the condition $\left\lVert \phi \right\rVert_2 = 1$. Originating in our case from the study of ground state energy asymptotics of the Pauli-Fierz Hamiltonian of quantum electrodynamics, problems of this type often emerge in different fields of mathematics such as optimization and image processing. By rescaling it is also closely related to the Nash inequality \cite{Nash} with minimizers corresponding to the saturation cases of the inequality. In fact Carlen and Loss \cite{CarlenLoss} give a proof of the inequality's sharpness using the Poincaré inequality and specific properties of eigenfunctions of the Neumann Laplacian on unit balls which are also shown to represent all cases of equality up to translation, rescaling and normalization. While the existence of minimizers in a slightly more general model has also been addressed in a recent work by Siegel and Tekin \cite{Siegel} using the Kolmogorov-Riesz theorem and Gagliardo-Nirenberg inequality, we would like to offer a more emphasized insight into the inherent issue of non-reflexivity and non-duality of the singular Lebesgue space $L^1$ at the heart of the problem. 

\begin{definition}
For any $\beta > 0$ let $\cF_\beta: H^1(\RR^3)\cap L^1(\RR^3) \to \RR_0^+$ be given by
\begin{align*}
\cF_\beta(\phi)  & := \frac{1}{2}\left\lVert\nabla\phi\right\rVert_2^2 + \beta \left\lVert \phi \right\rVert_1 = \int f_\beta(\partial_{x_\nu} \phi, \phi) \mathrm{d}x^3
\end{align*}
with $f_\beta\left(\partial_{x_\nu} \phi, \phi\right) := \frac{1}{2}\left\lvert \nabla \phi\right\rvert^2 + \beta \left\lvert \phi \right\rvert$
\end{definition}


\begin{lemma} \label{rearrangement}
Let $\phi^*$ be the spherically symmetric decreasing rearrangement of $\left\lvert\phi\right\rvert$. Then
\begin{align*}
\cF_\beta(\phi^*) \leq \cF_\beta(\phi) .
\end{align*}
\begin{proof}
The inequality follows from standard results on rearrangements such as $\left\lVert \phi^* \right\rVert_1 = \left\lVert \phi \right\rVert_1$ and the Pólya-Szeg\H{o} inequality $\left\lVert\nabla\phi^*\right\rVert_2 \leq \left\lVert\nabla\phi\right\rVert_2$
\end{proof}
\end{lemma}
%
%
%
%
Establishing existence of minimizers in the context of calculus of variations typically relies on a variety of arguments in a procedure known as the direct method. Its central component consists in the Banach-Alaoglu theorem or implications thereof proving sequential compactness of minimizing sequences in some weak-(*) topology, while a corresponding lower semicontinuity property guarantees that their limit attains the minimum.
In this case however the idiosyncrasies of the singular $p=1$ Lebesgue space impose additional technical difficulty onto the direct method. Not only is $L^1(\RR^3)$ not reflexive but also not a dual space which can for instance be seen by a combination of the Banach-Alaoglu and Krein-Milman theorems. The lack of a predual space renders the very notion of weak-* topology inapplicable and in terms of weak topology the Eberlein-\v{S}mulian theorem  even establishes that the closed unit ball of a Banach space $X$ is weakly sequentially compact if and only if $X$ is reflexive.

We present a solution to the problem of $L^1$-convergence by a introducing a series of adaptations to the method and argue as follows:
\begin{itemize}
\item For any finite $R$, the intersection $\left(H^1\cap L^1\right)(B_R(0))$ is a uniformly convex space when equipped with a suitable composite norm inducing an equivalent topology. 
\item Uniformly convex spaces are reflexive by virtue of the Milman-Pettis theorem \cite{Milman}.
\item The direct method applied to $\left(H^1\cap L^1\right)(B_R(0))$ gives rise to a weakly-* convergent subsequence with converging norms. In uniformly convex spaces this is equivalent to strong convergence (cf. \cite{Brezis} ch. 3).
\item Any minimizing sequence of $\cF$ is uniformly tight in $L^1$-norm.
\item Strong convergence on all finite balls and uniform tightness imply total boundedness. According to the general Bolzano-Weierstrass theorem completeness and total boundedness imply sequential relative compactness.
\end{itemize}

While this procedure allows to select a strongly $L^1$-convergent subsequence from any minimizing sequence, the usual direct method yields the analogous result in terms of $H^1$-topology. In conjunction one may conclude the existence of a minimizer.

Let us proceed to the first step. Following the works of Clarkson \cite{Clarkson} who introduced the concept of uniform convexity one can show that general product spaces are reflexive in the topology of a composite norm satisfying certain properties such as homogeneity and strict convexity as long as each factor is reflexive. The additional difficulty emerging here however, the non-reflexivity of $L^1$, can be dealt with on finite $B_R(0)$ (but not $\RR^3$) by relative bounds on $L^1$- in terms of $H^1$-norms. In this way the we can show that $(H^1\cap L^1)(B_R(0))$ is a uniformly convex Banach space with respect to a certain norm and therefore reflexive.

\begin{definition}
Let the norm $\left\lVert . \right\rVert_{H^1\cap L^1} : H^1(\RR^3)\cap L^1(\RR^3) \to \RR^+_0$ be defined as $\left\lVert \phi \right\rVert_{H^1\cap L^1}:= \sqrt {\left\lVert \phi \right\rVert_{H^1}^2 + \left\lVert \phi \right\rVert_{L^1}^2}$. Then $\left\lVert . \right\rVert_{H^1\cap L^1}$ is continuous, homogeneous, strictly increasing and strictly convex in $\left\lVert \phi \right\rVert_{H^1}$ and $\left\lVert \phi \right \rVert_{L^1}$. The last property signifies that
\begin{align*}
& \sqrt {(\left\lVert \phi \right\rVert_{H^1} +\left\lVert \eta \right\rVert_{H^1})^2 + (\left\lVert \phi \right\rVert_{L^1} +\left\lVert \eta \right\rVert_{L^1})^2} \\
&  < \sqrt {\left\lVert \phi \right\rVert_{H^1}^2 +\left\lVert \phi \right\rVert_{L^1}^2} + \sqrt {\left\lVert \eta \right\rVert_{H^1}^2 + \left\lVert \eta \right\rVert_{L^1}^2}
\end{align*}
for any $\phi, \eta$ unless $(\left\lVert \phi \right\rVert_{H^1},\left\lVert \phi \right\rVert_{L^1})=c\cdot(\left\lVert \eta \right\rVert_{H^1},\left\lVert \eta \right\rVert_{L^1})$ for some $c\geq 0$.
\end{definition}

\begin{definition}
A Banach space $\left(X,\left\lVert.\right\rVert\right)$ is said to be uniformly convex if
\begin{align*}
\forall \varepsilon >0 \; \exists \delta >0 \; & \forall x,y\in X, \left\lVert x\right\rVert = 1, \left\lVert y\right\rVert = 1: \\
\left\lVert x-y\right\rVert > \varepsilon \quad & \Rightarrow \quad \left\lVert \frac{x +y}{2}\right\rVert < 1-\delta
\end{align*}
The conditions $\left\lVert x\right\rVert, \left\lVert y\right\rVert = 1$ may be equivalently replaced with $\left\lVert x\right\rVert, \left\lVert y\right\rVert \leq 1$. 
\end{definition}

\begin{theorem}
Let $R<\infty$. Then $ \left(H^1\cap L^1(B_R(0)), \left\lVert . \right\rVert_{H^1\cap L^1}\right)$ is a uniformly convex, reflexive Banach space.
\begin{proof}
First of all, any Cauchy sequence with respect to $\left\lVert . \right\rVert_{H^1\cap L^1}$ is a Cauchy sequence in the Banach spaces $H^1$ and $L^1$. Since $H^1(\RR^3)\cap L^1(\RR^3)$ is dense in both $H^1$ and $L^1$ the two limits must coincide and there is a common limit in $H^1 \cap L^1$.
Now let $(\phi_i)_i, (\eta_i)_i \in H^1(\RR^3)\cap L^1(\RR^3)$  two sequences with $\left\lVert \phi_i \right\rVert_{H^1\cap L^1} = \left\lVert \eta_i \right\rVert_{H^1\cap L^1} = 1, \; \forall i$ and $\left\lVert \phi_i + \eta_i \right\rVert_{H^1\cap L^1} \to 2, \;i\to\infty$. We demonstrate that $\left\lVert \phi_i - \eta_i \right\rVert_{H^1\cap L^1} \to 0$ which implies uniform convexity. Now due to $\left\lVert . \right\rVert_{H^1\cap L^1}$ being strictly increasing 
\begin{align*}
2 \leftarrow \left\lVert \phi_i + \eta_i \right\rVert_{H^1\cap L^1} & = \sqrt {\left\lVert \phi_i + \eta_i \right\rVert_{H^1}^2 + \left\lVert \phi_i +\eta_i \right\rVert_{L^1}^2} \\
& \leq \sqrt {\left(\left\lVert \phi_i  \right\rVert_{H^1} + \left\lVert  \eta_i \right\rVert_{H^1}\right)^2 + \left(\left\lVert \phi_i \right\rVert_{L^1} + \left\lVert\eta_i \right\rVert_{L^1}\right)^2} \\ 
& \leq \sqrt {\left\lVert \phi_i  \right\rVert_{H^1}^2 + \left\lVert \phi_i \right\rVert_{L^1}^2} + \sqrt {\left\lVert  \eta_i \right\rVert_{H^1}^2 + \left\lVert \eta_i \right\rVert_{L^1}^2} \\
& = \left\lVert \phi_i \right\rVert_{H^1\cap L^1} + \left\lVert \eta_i \right\rVert_{H^1\cap L^1} = 2.
\end{align*}  
Since in both inequalities equality has to hold in the limit, strict convexity and continuity imply that 
\begin{align*}
\underset{i\to\infty}{\lim} \left(\left\lVert \phi_i  \right\rVert_{H^1} - \left\lVert \eta_i  \right\rVert_{H^1} \right) = 0, \quad \quad \underset{i\to\infty}{\lim} \left(\left\lVert \phi_i  \right\rVert_{L^1} - \left\lVert \eta_i  \right\rVert_{L^1} \right) = 0 .
\end{align*}
Assuming that $\left\lVert \phi_i - \eta_i \right\rVert_{H^1\cap L^1} \not\rightarrow 0$, by selecting a suitable subsequence we can infer without loss of generality that
\begin{align*}
\begin{split}
& \underset{i\to\infty}{\lim} \left\lVert \phi_i - \eta_i \right\rVert_{H^1\cap L^1} = \alpha > 0  \\
& \underset{i\to\infty}{\lim} \left\lVert \phi_i  \right\rVert_{H^1} = \underset{i\to\infty}{\lim} \left\lVert  \eta_i \right\rVert_{H^1} = \beta \\
& \underset{i\to\infty}{\lim} \left\lVert \phi_i \right\rVert_{L^1} = \underset{i\to\infty}{\lim} \left\lVert  \eta_i \right\rVert_{L^1} = \beta' \\
\end{split}
\begin{split}
& \underset{i\to\infty}{\lim} \left\lVert \phi_i - \eta_i \right\rVert_{H^1} = \gamma > 0 \\
& \left\lVert \phi_i  \right\rVert_{H^1} > 0, \quad \left\lVert \eta_i  \right\rVert_{H^1} > 0 \quad \forall i
\end{split}
\end{align*}
for some $\alpha, \beta, \beta', \gamma$. The penultimalte line is a consequence of the Hölder inequality on finite domains
\begin{align*}
\left\lVert f  \right\rVert_{L^1} = \int_{B_R(0)} |f| \mathrm{d}^3x \leq \mathrm{Vol}(B_R(0))^{\frac12} \cdot \left\lVert f  \right\rVert_{L^2} \leq \mathrm{Vol}(B_R(0))^{\frac12} \cdot \left\lVert f  \right\rVert_{H^1} .
\end{align*}
Next we define 
\begin{align*}
\widetilde{\phi}_i :=  \frac{\beta \phi_i}{\left\lVert \phi_i  \right\rVert_{H^1}}, \quad \widetilde{\eta}_i :=  \frac{\beta \eta_i}{\left\lVert \eta_i  \right\rVert_{H^1}}.
\end{align*}
Then $\underset{i\to\infty}{\lim} \left\lVert \phi_i - \widetilde{\phi}_i \right\rVert_{H^1} = \underset{i\to\infty}{\lim} \left\lVert \eta_i - \widetilde{\eta}_i \right\rVert_{H^1} = 0$ and 
$
\underset{i\to\infty}{\lim}  \left\lVert \widetilde{\phi}_i - \widetilde{\eta}_i \right\rVert_{H^1} = \gamma 
$.
But the space $H^1$ is uniformly convex itself, such that
\begin{align*}
\underset{i\to\infty}{\limsup}  \left\lVert \phi_i + \eta_i \right\rVert_{H^1} = \underset{i\to\infty}{\limsup}  \left\lVert \widetilde{\phi}_i + \widetilde{\eta}_i \right\rVert_{H^1} < 2\beta
\end{align*}
Since 
\begin{align*}
\sqrt{\beta^2 + \beta'^2 } = \sqrt{\left(\underset{i\to\infty}{\lim}  \left\lVert \phi_i \right\rVert_{H^1} \right)^2 + \left(\underset{i\to\infty}{\lim}  \left\lVert \phi_i \right\rVert_{L^1} \right)^2 } = \underset{i\to\infty}{\lim}  \left\lVert \phi_i \right\rVert_{H^1 \cap L^1} = 1
\end{align*} 
it follows that 
\begin{align*}
\underset{i\to\infty}{\limsup}  \left\lVert \phi_i + \eta_i \right\rVert_{H^1\cap L^1} & \leq \sqrt{\left(\underset{i\to\infty}{\limsup}  \left\lVert \phi_i + \eta_i \right\rVert_{H^1}\right)^2 + \left(\underset{i\to\infty}{\limsup}  \left\lVert \phi_i + \eta_i \right\rVert_{L^1}\right)^2 } \\
& < \sqrt{\left(2\beta\right)^2 + \left( 2\beta'\right)^2 } = 2\sqrt{\beta^2 + \beta'^2 } = 2 
\end{align*}
which is a contradiction.
It follows that $(H^1 \cap L^1)(B_R(0))$ is uniformly convex. Consequently it is reflexive by the Milman-Pettis theorem \cite{Milman} (cf. \cite{Brezis} ch. 3).
\end{proof}
\end{theorem}

\begin{lemma}
Let $(\phi_n)_n \subset \left(H^1\cap L^1\right)(\RR^3)$ with $\phi_n = \phi_n^*, \,\forall n\in\NN$ be bounded in $\left\lVert . \right\rVert_{H^1\cap L^1}$. Then $(\phi_n)_n$ is uniformly tight with respect to $\left\lVert . \right\rVert_{L^2}$, i.e.
\begin{align*}
\forall \varepsilon > 0 \; \exists R < \infty  \; \forall n\in\NN : \; \left\lVert \phi_n |_{\RR^3 \setminus B_R(0)} \right\rVert_{L^2} < \varepsilon .
\end{align*}
Additionally, if $(\phi_n)_n$ is a minimizing sequence of $\cF$, it is uniformly tight with respect to  $\left\lVert . \right\rVert_{L^1}$.
\begin{proof}
Since $\phi_n=\phi^*\in \left(H^1\cap L^1\right)(\RR^3)$ and the symmetric decreasing rearrangement is monotonic, $\phi_n$ is continuous. By a weak version of the fundamental theorem of calculus $\phi_n(r) = - \int_r^\infty \phi'_n(s) \mathrm{d}s $ for all $r\in \RR^+_0$. The $L^1$-integrability implies $\int_0^\infty \phi_n(r) r^2 \mathrm{d} r = - \int_0^\infty \int_r^\infty \phi'_n(s) r^2 \mathrm{d}s  \mathrm{d}r < \infty $ uniformly in $n$. \\
By Fubini's theorem however, the other iterated integral $\int_0^\infty \int_0^s \phi'_n(s) r^2 \mathrm{d}r \mathrm{d} s = \int_0^\infty \frac{s^3}{3} \phi'_n(s) \mathrm{d}s $ exists and is finite as well. In particular 
\begin{align*}
r^3 \phi_n(r) = - r^3 \int_r^\infty \phi'_n(s) \mathrm{d}s \leq - \int_r^\infty s^3\phi'_n(s) \mathrm{d}s \leq C < \infty 
\end{align*} 
uniformly in $n$ which shows $\phi_n(r)\leq \frac{C}{r^3}$ and therefore uniform tightness in $\left\lVert . \right\rVert_{L^2}$. 
Assuming that $(\phi_n)_n$ is not uniformly tight in $L^1$-norm, there is  some $\varepsilon > 0$ such that for all $R < \infty $ there exists some $n\in\NN$ with $\left\lVert \phi_n |_{\RR^3 \setminus B_R(0)} \right\rVert_1 \geq \varepsilon$. But since $(\phi_n)_n$ is uniformly tight in $L^2$-norm we can find a sequence $R_n \to \infty$ such that without loss of generality
\begin{align*}
\left\lVert \phi_n |_{\RR^3 \setminus B_{R_n}(0)} \right\rVert_1 & \geq \varepsilon & \mathrm{and} &   & \left\lVert \phi_n |_{\RR^3 \setminus B_{R_n}(0)} \right\rVert_2 < \frac{1}{n}, \quad \forall n\in\NN
\end{align*}
Next for any $n\in\NN$ we define $\widetilde{\phi}_n$ by 
\begin{align*}
\widetilde{\phi}_n(r) := \begin{cases} \phi_n(r) & r\leq R_n \\ (1-r)\phi_n(R_n) & R_n < r < R_n+1 \\ 0 & r \geq R_n +1 \end{cases}
\end{align*}
Since $\phi_n(r)\leq \frac{C}{r^3}$, by construction $\left\lVert \nabla \widetilde{\phi}_n |_{\RR^3 \setminus B_{R_n}(0)}\right\rVert_{2}^2 + \left\lVert \widetilde{\phi}_n|_{\RR^3 \setminus B_{R_n}(0)} \right\rVert_{1} \to 0$ and $\left\lVert \widetilde{\phi}_n|_{\RR^3 \setminus B_{R_n}(0)} \right\rVert_{2} \to 0$ as $n\to \infty$.
In conjunction with the above this implies
\begin{align*}
\underset{n\to\infty}{\lim} \left\lVert \widetilde{\phi}_n \right\rVert_{L^2}  & =  \underset{n\to\infty}{\lim} \left\lVert \phi_n \right\rVert_{L^2}  \\
  \underset{n\to\infty}{\lim} \left(\frac12 \left\lVert \nabla \widetilde{\phi}_n \right\rVert_{2}^2 + \beta \left\lVert \widetilde{\phi}_n \right\rVert_{1}\right) & < \underset{n\to\infty}{\lim} \left(\frac12 \left\lVert \nabla \phi_n \right\rVert_{2}^2 + \beta\left\lVert \phi_n \right\rVert_{1}\right) 
\end{align*}
which demonstrates that $(\phi_n)_n$ is not a minimizing sequence of $\cF$ and concludes the proof.
\end{proof}
\end{lemma}


We are now ready to show the existence of a minimizer in the spirit of the direct method of calculus of variations.

\begin{theorem}
There exists a minimizer $\phi\in H^1(\RR^3)\cap L^1(\RR^3)$ of $\cF_\beta$ under the condition $\left\lVert \phi \right\rVert_2 = 1$. The minimizer's spherically symmetric decreasing arrangement is unique and has compact support.
\begin{proof}
We consider the Banach space $H^1(\RR^3)\cap L^1(\RR^3)$ endowed with norm $\left\lVert \phi \right\rVert_{H^1\cap L^1}:= \sqrt {\left\lVert \phi \right\rVert_{H^1}^2 + \left\lVert \phi \right\rVert_{L^1}^2}$ where $\left\lVert \phi \right\rVert_{H^1}^2 = \left\lVert \phi \right\rVert_2^2 + \left\lVert \nabla\phi \right\rVert_2^2$. Under the condition $\left\lVert \phi \right\rVert_2 = 1$ the minimization of $\cF$ is equivalent to that of 
\begin{align*}
\widetilde{\cF}(\phi):= \frac{1}{2}\left(\left\lVert\nabla\phi\right\rVert_2^2 + \left\lVert \phi \right\rVert_2^2\right) + \beta \left\lVert \phi \right\rVert_1 = \frac{1}{2}\left\lVert \phi \right\rVert_{H^1}^2 + \beta \left\lVert \phi \right\rVert_{L^1}
\end{align*} 
$\widetilde{\cF}$ is bounded from below, so $\underset{\left\lVert\phi\right\rVert_2=1}{\inf} \, \widetilde{\cF}(\phi)$ exists. 
Let $(\phi_n)_n\in H^1(\RR^3)\cap L^1(\RR^3)$ be a minimizing sequence such that $\widetilde{\cF}(\phi_n) \to \underset{\left\lVert\phi\right\rVert_2=1}{\inf}\, \widetilde{\cF}(\phi)$. Transitioning to its symmetric decreasing rearrangement yields another smaller minimizing sequence and without loss of generality we may assume that $\phi_n =\phi_n^*$.

By definition, $\widetilde{\cF}$ is a smooth function of the $H^1$- and $L^1$-norms which are equivalent to the $H^1\cap L^1$- norm. Clearly $\widetilde{\cF}$ is coercive with respect to both $\left\lVert . \right\rVert_{H^1}$ and $\left\lVert . \right\rVert_{L^1}$, so $(\phi_n)_n$ is bounded in both norms. Additionally $\widetilde{\cF}$ is weak-* lower semincontinuous with respect to $H^1$.

An application of Banach-Alaoglu to the $H^1$-topology and the direct method procedure yield sequential weak-* compactness of $(\phi_n)_n$, weak-* convergence of a subsequence and convergence of norms all with respect to $\left\lVert . \right\rVert_{H^1}$. Since $H^1$ is a Hilbert space this automatically implies strong convergence of a subsequence and there is $\phi \in H^1(\RR^3)$ such that without loss of generality $\left\lVert \widetilde{\phi} -\phi_n\right\rVert_{H^1} \to 0$ and consequently also $\left\lVert \widetilde{\phi} -\phi_n\right\rVert_{L^2} \to 0$. In remains to prove the assertion in terms of $L^1$-convergence. 

For any $R < \infty$ the space $((H^1\cap L^1)(B_R(0)), \left\lVert . \right\rVert_{H^1\cap L^1})$ is uniformly convex {\&} reflexive and we consider the restrictions $\phi_n\vert_{B_R(0)}$. By the Banach-Alaoglu theorem, any bounded closed ball in $\left(H^1\cap L^1\right)(B_R(0))$ is sequentially weakly-* compact so there exist $\widetilde{\phi}\vert_{B_R(0)} \in \left(H^1\cap L^1\right)(B_R(0))$ such that $\phi_n\vert_{B_R(0)} \overset{\mathrm{w}-*}{\to} \widetilde{\phi}\vert_{B_R(0)}$. 
Moreover $\widetilde{\cF}(\phi)$ is weak-* lower semicontinuous so that $\widetilde{\cF}(\phi_n\vert_{B_R(0)}) \to \widetilde{F}(\phi\vert_{B_R(0)})$. Since $\phi \mapsto \widetilde{\cF}(\phi)= \frac{1}{2}\left\lVert \phi \right\rVert_{H^1}^2 + \beta \left\lVert \phi \right\rVert_{L^1} $ is strictly convex on the space of spherically symmetric decreasing $\phi = \phi^*$, it follows that $\left\lVert \phi_n\vert_{B_R(0)}\right\rVert_{H^1} \to \left\lVert \phi\vert_{B_R(0)} \right\rVert_{H^1}$ and $\left\lVert \phi_n\vert_{B_R(0)}\right\rVert_{L^1} \to \left\lVert \phi\vert_{B_R(0)} \right\rVert_{L^1}$.
As $(H^1\cap L^1)(B_R(0))$ is uniformly convex and reflexive, weak-* convergence and convergence of norms imply strong convergence  $\phi_n\vert_{B_R(0)} \overset{H^1\cap L^1}{\to}   \widetilde{\phi}\vert_{B_R(0)}$ and therefore $\phi_n\vert_{B_R(0)} \overset{L^1}{\to}   \widetilde{\phi}\vert_{B_R(0)}$. 

$L^1(\RR^3)$ is a Banach space and in particular complete. To prove that the minimizing family  $(\phi_n)_n\subset H^1\cap L^1(\RR^3)$ is relatively compact in $L^1$-norm it suffices to show that it is totally bounded, i.e. admits a finite covering of $\varepsilon$-balls for any $\varepsilon > 0$. But this follows from strong convergence $\left\lVert (\phi - \phi_n)|_{B_R(0)} \right\rVert_{L^1} \to 0$ on $H^1\cap L^1(B_R(0))$ and the fact that the minimizing family $(\phi_n)_n\subset H^1\cap L^1(\RR^3)$ is uniformly tight in $L^1$-norm. Hence there exists a $L^1$-convergent subsequence of $(\phi_n)_n\subset H^1\cap L^1(\RR^3)$. As $L^p$-convergence implies pointwise convergence a.e. the $H^1$- and $L^1$-limits coincide and by selecting subsequences
$
\left\lVert \widetilde{\phi} -\phi_n\right\rVert_{H^1} \to 0$ and $\left\lVert \widetilde{\phi} -\phi_n\right\rVert_{L^1} \to 0
$.
We conclude that there is a $\phi \in H^1\cap L^1(\RR^3)$ such that $\left\lVert \phi - \phi_n \right\rVert_{H^1\cap L^1(\RR^3)} \to 0$ with $\widetilde{\cF}(\phi) = \underset{\phi\in H^1\cap L^1,\;\left\lVert\phi\right\rVert_2=1}{\inf} \widetilde{\cF}(\phi) $. Moreover strong $H^1$-convergence implies strong $L^2$-convergence and therefore $\left\lVert \phi \right\rVert_{L^2} = \lim_{n\to\infty} \left\lVert \phi_n \right\rVert_{L^2} = 1$
\end{proof} 
\end{theorem}

\begin{theorem}
The unique spherically decreasingly rearranged and normed minimizer $\phi$ of $\cF_\beta$ is a solution to the Helmholtz equation on $B_R(0)$ with Neumann (and Dirichlet) boundary conditions
\begin{align*}
\begin{split} \left. \begin{cases} \left(\Delta +  \mu^2 \right) \phi(x) = \beta  , & |x| < R \\ \phi(x) = 0 , & |x|\geq R \end{cases} \right\rbrace  \end{split}, & \partial_{\vec{n}} \phi(x) = 0, \; x \in \partial B_R(0)
\end{align*}
for some parameters $\mu, R \in  \RR^+$
\begin{proof}
The symmetric decreasing rearrangement $\phi^*$ is always lower semicontinuous and we define $R$ by $R:=\inf \lbrace |x|  \mid   \phi^*(x)=0 \rbrace \in \RR^+$ or by $R:=\infty$ if the former set is empty. $\phi^* = \phi^*(|x|)$ is monotonic and therefore continuous up to at most countably many jump discontinuities. But since $\phi^*\in H^1(\RR^3)$ it follows that that $\phi^*$ is continuous. 
For any arbitrary but fixed $R$ we consider the set of admissible support-preserving variations in the point $\phi$ which is given by $\cD  :=  C_c^\infty \left( B_R(0); \RR \right)$.
If $\phi$ is a minimizer then there is a Lagrange multiplier $\lambda$ such that the total variation (with fixed $R$)  
\begin{align*}
\delta(\cF+\lambda \cG) = \left(\frac{\partial}{\partial \phi}(\cF+\lambda \cG)\right)(\delta\phi) + \sum_\nu \left(\frac{\partial}{\partial (\partial_{x_\nu}\phi)}(\cF+\lambda \cG)\right)(\delta\partial_{x_\nu}\phi)
\end{align*}
satisfies
\begin{align*}
0 & =\delta (\cF+\lambda \cG) = \left\langle \nabla \phi , \nabla (\delta \phi) \right\rangle + \beta \left\langle \one_{\phi > 0}, \delta \phi \right\rangle + \lambda\left\langle \phi, \delta \phi \right\rangle 
\end{align*}
for any $\delta \phi \in \cD$ in the sense of a Fréchet derivative. 
This gives rise to the Euler-Lagrange equation 
\begin{align*}
\frac{\partial (\cF+\lambda \cG)}{\partial \phi } = \sum_\nu \frac{\partial}{\partial x_\nu} \frac{\partial (\cF+\lambda \cG)}{\partial (\partial_{x_\nu} \phi)} 
\end{align*}
By the du Bois-Reymond theorem of calculus of variations (\cite{Gelfand} p 17f.) since $\phi \in H^1(\RR^3)$ and $\left(\frac{\partial^2}{\partial_{x_\nu}\phi \partial_{x_\mu}\phi} (f+\lambda g)\right)_{\nu, \mu}= \left(\frac{\partial^2}{\partial_{x_\nu}\phi \partial_{x_\mu}\phi} f\right)_{\nu,\mu} = 2\cdot \one > 0$ is positive definite everywhere, it follows that any solution to the Euler-Lagrange equation and therefore the minimizer $\phi$ has to be twice (weakly) differentiable.
Therefore for all $\delta \phi \in \cD$
\begin{align*}
0 & =\delta (\cF+\lambda \cG) = \left\langle \nabla \phi , \nabla (\delta \phi) \right\rangle + \beta \left\langle \one_{\phi > 0}, \delta \phi \right\rangle + \lambda\left\langle \phi, \delta \phi \right\rangle \\
& = \left\langle -\Delta \phi , \delta \phi \right\rangle + \beta \left\langle \one_{\phi > 0}, \delta \phi \right\rangle + \lambda \left\langle \phi, \delta \phi \right\rangle
\end{align*}
which implies the more explicit Euler-Lagrange equation
\begin{align*}
(- \Delta + \lambda) \phi + \beta = 0
\end{align*}
weakly on $B_R(0)$. 
Now $(- \Delta + \lambda) \phi + \beta = 0$, $\phi\in H^1(\RR^3)$ and existence of the second weak derivative imply that $\phi|_\Omega\in H^2(\Omega)$ for all bounded measurable domains $\Omega$. 
As a consequence of Gauss' divergence theorem in Sobolev spaces it follows
\begin{align*}
 \beta \cdot \mathrm{Vol}\left(B_{R'}(0)\right) & =  \int_{B_{R'}(0)} \Delta \phi - \lambda \phi \; \mathrm{d}^3x \\ 
 & = \int_{\partial B_{R'}(0)} \vec{\nabla}\phi \cdot \vec{n}\mathrm{d}S + \int_{B_{R'}(0)}  - \lambda \phi \; \mathrm{d}^3x
\end{align*}
for all $R' < \infty$. But since $\phi=\phi^*$, $\vec{\nabla}\phi \cdot \vec{n}$ and the surface integral are non-positive. With $\beta>0$ and taking the limit $R' \nearrow R $ it follows that $ \beta \cdot \mathrm{Vol}\left(B_R(0)\right) \leq |\lambda| \cdot \left\lVert \phi \right\rVert_1 < \infty$ hence $R < \infty $. In particular $\phi \in H^2(\RR^3)$.
Finally $\lambda < 0$ has to hold, for if $\lambda \geq 0$, then $\Delta \phi = \lambda \phi + \beta > 0$ and $\phi$ would be weakly subharmonic. By the maximum principle (in Sobolev space, c.f. \cite{Brezis} ch. 9) $\phi | _{\partial B_R(0)} = 0$ would imply $\phi \leq 0$ which is a contradiction. Alternatively this also follows from the above divergence theorem.

Finally it follows from a variation of $R$ in the derivation of the Euler-Lagrange equation that $\partial_r \phi (R) = 0$ since in the sense of a Fréchet derivative
\begin{align*}
0  = \frac{\delta}{\delta R} (\cF+\lambda \cG)|_{\phi=\phi^*} & = \left\lvert \partial B_R(0)\right\rvert \cdot \left( \left(\nabla \phi (R)\right)^2 + \beta \phi(R)  + \lambda \phi(R)^2 \right) \\
& = 4\pi R^2 \left(\partial_r \phi (R) \right)^2
\end{align*}
\end{proof}
\end{theorem}

\begin{theorem}
The unique spherically decreasingly rearranged and normed minimizer $\phi$ of $\cF_\beta$ with is given by
\begin{align*}
\phi(r) = \begin{cases} a \cdot\frac{\sin( \mu r)}{\mu r} + \frac{\beta}{\mu^2}, & r \leq R \\ 0, & r>R \end{cases}
\end{align*}
for some parameters $a,\mu, R \in \RR^+$. The parameters are implicitly uniquely determined by $R = \min \lbrace r > 0 \mid \phi(r)=0 \rbrace$, $\phi'(R)=0$ and $\left\lVert \phi \right\rVert_2 = 1$.
\begin{proof}
Let $\mu^2 := - \lambda$ as above. The partial differential equation $(- \Delta +  \mu^2) \phi + \beta = 0$ with $\phi>0$ on $B_R(0)$ that constitutes the Euler-Lagrange equation is the Helmholtz equation in dimension 3. It admits a characterization of homogeneous solutions in terms of spherical Bessel functions $j_l,y_l$ and spherical harmonics $Y_l^m$ according to
$
\phi_\mathrm{hom} = \sum_{l=0}^\infty \sum_{m=-l}^{l} \left(a_{lm}j_l(\mu r) + b_{lm}y_l(\mu r)\right)Y_l^m(\vartheta, \varphi) .
$
The spherically symmetry of $\phi=\phi^*$ implies $a_{lm}=0=b_{lm}$ for any $(l,m)\neq(0,0)$, while the inhomogeneous solution is easily seen to be given by $\phi_\mathrm{inhom} =  \frac{\beta}{\mu^2}$. It follows that 
$
\phi =  a_{00} j_0(\mu r) + b_{00} y_0(\mu r) + \frac{\beta}{\mu^2}.
$
But the regularity of $\phi$ in $r=0$ requires $b_{00}$ to vanish since $\underset{r\to 0}{\lim}\, y_0(r) = -\infty$ for any Bessel functions of the second kind. Thus
\begin{align*}
\phi(r) =  a_{00}j_0(\mu r) + \frac{\beta}{\mu^2}= a_{00} \cdot\frac{\sin( \mu r)}{\mu r} + \frac{\beta}{\mu^2}.
\end{align*} 
for some $a_{00} =:a \in \RR^+$ such that $R$ is the smallest zero of $\phi$. It is then easy to verify that the parameters $a,\mu,R \in \RR^+$ are uniquely determined by the three conditions $R = \min \lbrace r > 0 \mid \phi(r)=0 \rbrace$, $\phi'(R)=0$ and $\left\lVert \phi \right\rVert_2 = 1$ 
\end{proof}
\end{theorem}

The following virial theorem is helpful in determining minimizing parameters $a,\mu,R$ for an arbitrary $\beta > 0$

\begin{lemma} \label{virial}
Let $\phi$ be the minimizer of $\cF_\beta$. Then 
$
 \left\lVert\nabla\phi\right\rVert_2^2 = \frac32 \beta \left\lVert \phi \right\rVert_1
$
\begin{proof}
Let $\nu \in \RR^+ $ be a relative scaling factor and $\phi_{(\nu)}(x):=\nu^{-\frac32}\phi(\nu^{-1} x)$ be a $L^2$-unitary rescaling. Since both $\left\lVert\nabla\phi\right\rVert_2^2$ and $\left\lVert \phi \right\rVert_1$ are homogeneous functionals and smooth under rescaling and by the minimality of $\phi$ 
\begin{align*}
0 & = \frac{\mathrm{d}}{\mathrm{d}\nu} \cF(\phi_{(\nu)}) \Bigm\vert_{\nu=1} = \frac{\mathrm{d}}{\mathrm{d}\nu} \left(\frac{1}{2\nu^2}\left\lVert\nabla\phi\right\rVert_2^2 +  \beta \nu^{\frac32} \left\lVert \phi \right\rVert_1 \right)\Bigm\vert_{\nu=1} \\
& = - \left\lVert\nabla\phi\right\rVert_2^2 + \frac32\beta  \left\lVert \phi \right\rVert_1 .
\end{align*}
\end{proof}
\end{lemma}

For any $\beta >0 $ we denote by $\alpha_\beta, \mu_\beta, R_\beta$ the parameters pertaining to the minimizer of $\cF_\beta$.

\begin{proposition} \label{parameters}
There exist constants $a_1, \mu_1, R_1 \in \RR^+ $ such that for any $\beta>0$ 
\begin{align*}
a_\beta = \beta^\frac37 a_1, \quad \mu_\beta = \beta^\frac27 \mu_1, \quad R_\beta = \beta^{-\frac27} R_1.
\end{align*}
Moreover, there is a constant $F_1 \in \RR^+ $ such that $F_\beta = \beta^{\frac47}F_1$
\begin{proof}
Let $\nu \in \RR^+$ be a relative scaling factor and $\phi_{(\nu)}(x):=\nu^{-\frac32}\phi(\nu^{-1} x)$. Then for any $\beta >0$
\begin{align*}
\cF_\beta(\phi_{(\nu)}) = \frac{1}{2\nu^2}\left\lVert\nabla\phi\right\rVert_2^2 +  \beta \nu^{\frac32} \left\lVert \phi \right\rVert_1 = \frac{1}{\nu^2}\left(\frac12\left\lVert\nabla\phi\right\rVert_2^2 +  \beta \nu^{\frac72} \left\lVert \phi \right\rVert_1 \right) = \frac{1}{\nu^2} \cF_{\beta \nu^{\frac72}} (\phi) .
\end{align*}
Therefore we have the equivalence
\begin{align*}
\phi_{(\nu)} \text{ is the minimizer of } \cF_\beta \quad & \Leftrightarrow \quad \phi \text{ is the minimizer of } \cF_{\beta \nu^{\frac72}}
\end{align*}
Now if $\beta, \beta' > 0$, then choosing $\nu >0$ such that $\beta' = \beta \nu^{\frac72}$ implies 
\begin{align*}
R_{\beta'} = R(\phi) =  \nu^{-1} R(\phi_{(\nu)})  = \nu^{-1} R_\beta 
\end{align*}
hence $ \left(\frac{R_\beta}{R_{\beta'}}\right)^{\frac72} = \nu^\frac72 = \frac{\beta'}{\beta} $.
Furthermore we have for any $\beta, \beta'>0$
\begin{align*}
 a_{\beta'}+\frac{\beta'}{\mu_{\beta'}^2}  = \phi(0) & = \nu^\frac32 \phi_{(\nu)}(0) =  \nu^\frac32 \left(a_{\beta}+\frac{\beta}{\mu_{\beta}^2}\right) \\
  a_{\beta'}\frac{\sin( \mu_{\beta'} R_{\beta'})}{\mu_{\beta'} R_{\beta'}} +\frac{\beta'}{\mu_{\beta'}^2} = 0 & = \nu^\frac32\cdot 0 =  \nu^\frac32 \left( a_{\beta}\frac{\sin( \mu_{\beta} R_{\beta})}{\mu_{\beta} R_{\beta}} +\frac{\beta}{\mu_{\beta}^2}\right)
\end{align*}
which is only solvable if $\frac{\sin( \mu_{\beta} R_{\beta})}{\mu_{\beta} R_{\beta}}= \mathrm{sinc} (\mu_{\beta} R_{\beta})$ is constant, i.e.
$
\mu_{\beta'} R_{\beta'} \equiv \mu_{\beta} R_{\beta} = \text{const.}
$
With $\nu^\frac72 = \frac{\beta'}{\beta}$ this implies
$\frac{\mu_{\beta'}}{\mu_{\beta}} = \nu$ and $\frac{a_{\beta'}}{ a_{\beta}} = \nu^\frac32$
so that in summary there exist constants $a_1, \mu_1, R_1 > 0$ such that for any $\beta>0$
\begin{align*}
a_\beta = \beta^\frac37 a_1, \quad \mu_\beta = \beta^\frac27 \mu_1, \quad R_\beta = \beta^{-\frac27} R_1.
\end{align*}
The second assertion then follows from Lem. \ref{virial} and noting that $\left\lVert\nabla\phi\right\rVert_2^2 \sim R_\beta^{-2}$. The fact that $F_1>0$ can be seen by computation and $a_1, \mu_1, R_1 >0$ or directly derived from the Nash inequality \cite{Nash}, stating that $\left\lVert \phi \right\rVert_2^{2+4/n} \leq  C_n\left\lVert \nabla \phi \right\rVert_2^2 \left\lVert \phi \right\rVert_1^{4/n}$ for some only dimension dependent constant $C_n > 0$, and a weighted arithmetic geometric mean inequality.
\end{proof}
\end{proposition}


\inputencoding{latin2}
\bibliographystyle{plain}
\inputencoding{utf8}

\end{document}